# *Timaeus*' Puzzle of the Innumerable Interstices in a Universe Without Void

Luc Brisson
Centre Jean Pépin
CNRS-ENS/UMR 8230

Salomon Ofman
Institut mathématique de Jussieu-
Paris Rive Gauche/HSM
Sorbonne Université-Université de Paris

**Abstract.** Some of the most challenging problems in *Timaeus*' cosmology arise from the geometry of a universe without any void. On the one hand, the universe is spherical in shape; on the other hand, it must be entirely filled with the four basic particles that make up all bodies in the universe, each shaped like one of four regular polyhedra (cubes, tetrahedra, octahedra and icosahedra). The faces of all these particles are composed of right triangles. However, this leads to two mathematical impossibilities.

1. Obtaining a spherical surface from linear surfaces, as it is impossible to create a circle from straight lines.
2. Obtaining a complete tiling of a sphere using regular polyhedra, without any voids or intersections between these polyhedra.

The first problem is addressed in another article slated for publication. In the present one, our focus will be on the second problem, for which we will present a solution within the framework of Timaeus' cosmology. The crux of this solution lies in a feature of *Timaeus*' universe that sets it apart from almost all ancient cosmologies. Instead of being composed of rigid parts, it is a dynamic living body in which all basic components are in constant motion, continuously undergoing both destruction and reconstruction.

In the first part, we examine the main features of *Timaeus*' cosmology relevant to our issue. In the second part, we analyze the paradox in details and its consequences for *Timaeus*' cosmology. We then discuss the common 'solutions' and highlight their shortcomings. Finally, we propose a solution that we believe is consistent with Plato's text and independent of the choice of the major schools of interpretation of the *Timaeus*.

## I. Introduction
### 1) A brief summary of the beginning of Timaeus' account.

At the beginning of his account, Timaeus states that the ordered universe requires the intervention of the *demiourgos*,[1] who must contend with 'necessity' ('ἀνάγκη', 47e-48d), from which he encounters some resistance. Because the god is inherently good and devoid of jealousy (29e), the *demiourgos* creates the universe as the most perfect and beautiful entity possible. However, due to the resistance encountered, the *demiourgos* is not all-powerful, and the realization depends on 'necessity', meaning that it is achieved to the extent that 'necessity' allows, either 'voluntary or under persuasion'.[2] Since the primary objective of the *demiourgos* is to create the 'most beautiful and perfect' universe, the outcome is always the 'most beautiful and perfect' *possible*, not an absolute perfection. The result is an ordered totality, which is a living creature with a soul and a body (30b). To be the most beautiful and perfect *possible*, the universe needs a model, which is the living creature of the absolute eternal immutable intelligible Forms (30c). However, the universe is everlasting not by its nature, but by the will of the *demiourgos*' (32c2-4; 37d3-4; 41a8-b6). It moves according to the motion most akin to immobility, which is rotation along an axis, and it possesses the most beautiful shape possible: the spherical one. More generally, the limitation 'insofar as is possible', which imposes constraints on what the *demiourgos* can produce, since there are impossible things even for the god, is an essential feature for understanding of *Timaeus*' cosmology.

---

[1] One meaning of *demiourgos* is precisely 'craftsman'.
[2] 'ὅπηπερ ἡ τῆς ἀνάγκης ἑκοῦσα πεισθεῖσά τε φύσις ὑπεῖκεν' (56c5-6); also 47e-48a.

The passages we summarized above have been the subject of extensive debate among Platonist commentators, from Antiquity to the present day. Scholars are deeply divided on whether to interpret these passages as metaphorical, literal or mythical. [3] The meaning of the *demiourgos*, the soul of the universe, and the role of 'necessity' in opposition to the *demiourgos* are also highly controversial. All these fundamental questions are beyond the scope of this article. In fact, the issues addressed here do not depend on choosing one standpoint over another. Consequently, the interpretation and solution presented here does not depend on the various fundamental understandings of the *Timaeus*, although they may be helpful for supporting (or criticizing) one argument or another.

## 2) The paradoxes

The focus of this article is to examine two aspects of *Timaeus*' that seem to result in a geometrical contradiction. On the one hand, the *demiourgos* gives the universe the most perfect shape — a spherical form (33b, 62e-63a). On the other hand, he emphasizes that no void space ('κενὴν χώραν', 58a7) or simply no void ('κενόν') exists [4] within or outside it (60c1; 79b1, c1; 80c3). However, Timaeus asserts that all bodies in the universe are composed of four basic particles — fire, air, water and earth — each shaped as one of the four regular polyhedra: tetrahedra, octahedra, icosahedra, and cubes, respectively.

- Aristotle already noted that it is impossible to fill a space, [5] let alone a sphere, [6] with these four polyhedra without leaving voids between them (*On the Heavens* VIII, 1, 306b6). He also pointed out that Plato's cosmology was flawed, as it contradicts mathematics, 'the most exact science' (*Ibid.*, III, 7, 306a28)
- Moreover, while asserting several times that no void can exist in the universe, Timaeus also claims that there are numerous 'interstices' ('διάκενα') between the basic particles.

Most commentators, following Aristotle, have considered Plato's cosmology to be contradictory and even not to be taken seriously. [7] In this article, we aim to demonstrate that, due to the unique features of *Timaeus*' cosmology, which differ significantly from Aristotle's, both the mathematical contradiction and the seemingly internal inconsistency can be resolved.

---

[3] For the first interpretation, see, for instance Ferrari (2022), p. xlii-xliii; for the second, Broadie (2012); for the third, we refer to Burnyeat (2005).

[4] 'there is no such thing as a void' ('τὸ κενὸν μηδὲν εἶναι', 79c1).

[5] Aristotle does not define what he means by 'space' ('τόπος'). However, it seems likely that the best interpretation is the simplest one: a portion of the total space whose boundary is continuous — more precisely, connected in the modern sense, meaning it forms a single whole. Since Aristotle asserts that this is only possible for cubes and tetrahedra, he implies the impossibility for octahedra and dodecahedra tiling *any* space, whatsoever (in the aforementioned sense). For instance, cubes can only tile parallelepipeds or combination of parallelepipeds, not general spaces. In particular, they cannot tile the sphere. Aristotle makes an error when he claims that tetrahedra can fill a space (at least according to the interpretation of 'space provided above).

[6] The term 'sphere' in ancient Greek, as well as sometimes today, is ambiguous. It can refer either to the points inside the sphere (what we would call a 'ball' in modern terms) or to its boundary (what we now call a 'sphere'). Since the meaning is clear from the context here, we will use 'sphere' as it was used in ancient Greece. The same ambiguity exists for the term 'circle', which can refer either to the points inside the circumference or to the circumference itself. For the same reason, we will maintain this ambiguity as well.

[7] For a review on this matter, see LloydA (1968), p. 78-79.

The key point is that *Timaeus*' universe, as a living being, is not rigid but rather has all its parts continuously moving and transforming.

## II. *Timaeus*' universe
### 1. Shaping the universe

Timaeus' universe is a living being ('ζῷον'), with both a body and a soul (36e). However, since it encompasses the totality of all living beings, it differs from them in various ways. First, because there is nothing external to the universes, it does not rely on anything beyond itself; it is self-sufficient. This is why it does not need senses, organs or limbs (33b-34a), allowing its body to be shaped in the most symmetrical form, [8] a sphere (33b). Consequently, the *demiourgos* gave the universe a spherical shape, and its soul the shape of two concentric circles with the same radius but tilted with respect to each other. One of these circles is called the 'circle of the Same', and the other, the 'circle of the Different' (36b *ff*). The *demiourgos* then extends the soul to cover the entire body of the universe, from its outermost edge to its center and beyond (34b3-4). In the remainder of the text we will use the term 'universe' to refer to the body of this 'living' being; otherwise, we will explicitly refer to the 'soul of the universe'.

### 2. The constitution of the universe

Timaeus asserts that the universe is completely filled with minuscule particles, too small to be seen (56c1-2), leaving 'no empty *space* to be left over.' [9] In essence, there is absolutely no void in the universe; it is entirely occupied by material bodies. These bodies are composed of four types of basic particles, each having the shape of one of the four regular polyhedra: tetrahedra, octahedra, icosahedra, and cubes (56c-57c). The particles are formed by assembling their faces; four, eight and twenty equilateral triangles for the first three polyhedra, respectively, and six squares for the cubes. The equilateral triangles are composed of six half-equilateral triangles and the squares are composed of four right isosceles triangles (53b-55d).

All the bodies in the universe are composed of these four types of basic particles distinguished by their geometrical forms. As noted in I.2) *supra*, they are shaped as four regular polyhedra: the tetrahedron corresponds to fire, the octahedron to air, the icosahedron to water, and the cube to earth. These *material* particles are formed from *mathematical* triangles, with the faces of the former three particles being formed by six half equilateral triangles, and the faces of the latter by four isosceles right triangles.

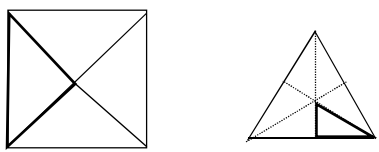

Figure 1

These basic particles are distinguished by the type of the basic triangles composing their faces. The particles of earth differ from the other three types because their faces of are formed by isosceles triangles, whereas the faces of the others are formed by half-equilateral triangles. These remaining particles are further distinguished by the number of basic half-equilateral triangles forming their faces: 24 for the tetrahedra, 48 for the octahedra and 120 for the dodecahedra. Thus, the particles are characterized by their faces, which are their boundaries, rather than their volumes. Additionally, they are in constant motion, and due to the absence of void, they continuously disintegrate into and reconstitute from the basic triangles (the isosceles or half-equilateral right triangles) (56c-57c). [10]

---

[8] The term 'σύμμετρος' means 'commensurable' in ancient Greek mathematics, but more generally can also mean 'in due proportion', and 'symmetrical'; See, for example, LloydB (2010), p. 455.
[9] 'κενὴν χώραν οὐδεμίαν ἐᾷ λείπεσθαι.' (58a7).

## 3. The spherical form of the universe

The ordered universe of the *Timaeus* requires the intervention of the *demiourgos* (30a-c) who always seeks to produce the most beautiful and perfect things. However, the *demiourgos* must contend with 'necessity' ('ἀνάγκη', 47e-48d), which presents some resistance. In any case, an image cannot achieve the same perfection as its model (37d3-4). This limitation appears so ingrained in Timaeus' construction that the *demiourgos* often does not even mention it, instead referring simply to the perfect universe and its perfect shape, as in the passage 32d-34b. Therefore, the mathematical paradox is easily resolved within Timaeus' framework. While the universe is spherical in shape, it is not an exact mathematical sphere but is crafted to be as *similar* to one as possible (34b1; 29a7). [11]

## 4. The movement of the universe

In fact, the universe itself is in motion. However, since nothing external to the universe exists, not even space, [12] as it encompasses the totality of physical space, its movement cannot involve *any* kind of translation. Since the universe cannot be as perfect as its intelligible model, it cannot be absolutely at rest. Therefore, the *demiourgos* gives it the most beautiful and perfect motion possible, one closest to reason and intelligence: rotation upon itself (34a; 40b). Moreover, since there is no space outside the universe, its rotation cannot result in any change of location. Consequently, the universe cannot have an arbitrary shape, such as being a polyhedron, as Aristotle noted. [13] This claim is presented as common knowledge in geometry, and was certainly well-known during Plato's time. However, Aristotle's assertion that only the sphere has the property to move without changing location is incorrect. In fact, this holds true for any solid of revolution, such as a cone or a cylinder. [14] While Plato indeed asserts that the universe rotates upon itself without any change place (34a1-4, 37c6), he does not make Aristotle's incorrect claim. The *demiourgos* gives the universe its spherical shape because the sphere is the most beautiful and perfect solid, not because it is only one capable of moving without changing location. In the *Laws*, when the Athenian painstakingly describes all kinds of motion, the second type is defined as rotation along an axis, with the axis itself also moving, corresponding to the motion of a wheel on a moving chariot. [15] However, a wheel can rotate on itself without changing location, whether it is not mounted on a chariot or the chariot is at rest for repairs. In fact, the wheel is a solid of revolution, essentially a cylinder. Thus, it is likely that Plato was aware that solids, other than spheres, can move without changing their location.

---

[10] Timaeus explains numerically how particles of fire, air, and water can be transformed into each other, while particles of earth cannot undergo such transformation because their faces are made of different triangles (56c-57c). For details on Timaeus' unusual construction, see, for instance, Brisson-Ofman (2021). For a difference approach, see Gregory (2008), in particular p. liii.

[11] While this is not the primary focus here, we will address this question and the mathematical difficulties it leads to in a forthcoming article.

[12] In the usual ancient Greek meaning, often synonymous with place, the 'common sense concept of space (space as a three-dimensional extension)', as argued by Algra (Algra (1995), p. 102); see also Johansen (2004), p. 188 for its different uses in the *Timaeus*.

[13] *On the Heavens* II, 4, 287a14-22.

[14] As Sorabji notes, Aristotle curiously excludes the shape of an egg (*On the Heavens* II, 4, 287a21-23), despite its ability to rotate around its vertical axis without changing its location (Sorabji (1988), p. 135). This exclusion makes it difficult to dismiss Aristotle's statement as merely a slip of the pen.

[15] For the first motion (without changing place): '- *And some will do this in one location, and others in several.* - *You mean, we will say, that those which have the quality of being at rest at the center move in one location, as when the circumference of circles that are said to stand still revolves?*' (cf. X, 893c4-7); for the second one (the axis also moving): '*you seem to me to mean all things that move by locomotion, continually passing from one spot*

Contrary to Aristotle's assertion of the necessity for the universe to be a sphere, Timaeus presents it as the deliberate choice of the *demiourgos*, since the god always seeks to produce the most beautiful and perfect things. The choice of a sphere rotating on its own axis arises from its status as the most beautiful and perfect solid in motion. However, as we have emphasized, the *Timaeus* is not purely a philosophical fiction; it also addresses these issues from an astronomical standpoint. In particular, mathematical principles apply and, in some sense, limit the power of the *demiourgos*.

The fact that the universe in its totality to be in motion is an important feature that distinguishes Timaeus' cosmology from other models, including Aristotle's. In most models of the universe, the motion of the different bodies within it is generally mutually independent of each other, but in Timaeus' cosmology, the movement of the universe is truly 'universal', meaning that all its components rotate together (34a; 36c). Thus, the usual view of a static space, within which some bodies move while others remain at rest, must be replaced by a dynamic model, where the universe in motion carries all its components with it.

### III. The dynamics of the universe
### 1. The war between the particles

The way that the particles are formed from and destroyed into the two basic right triangles is explained in 56d1-e7. The striking feature is that Timaeus is not interested in the question of volumes (cf. II.2, *supra*). To the astonishment of many modern commentators, the 'laws of transformation' between particles — for instance, water into air and vice-versa — depend exclusively of the number of the two basic right triangles composing them. Only this number is considered, while their volumes are simply disregarded and may vary widely. Another striking feature of the passage is its warlike overtones, with the particles of each of the four primary bodies battling against the others. [16] Timaeus describes in detail the encounter of particles with those of fire and of air or water:

> *whenever a small amount of fire is enveloped ('περιλαμβανόμενον') by a large quantity of air or water or perhaps earth and is agitated inside them as they move ('φερομένοις'), and is defeated ('νικηθέν') in its struggle ('μαχόμενον') and shattered to bits ('καταθραυσθῇ'), then any two fire corpuscles may combine to constitute a single form of air. And when air is overpowered ('κρατηθέντος') and broken down ('κερματισθέντος'), then two and one half entire forms of air will be consolidated into a single, entire form of water.* (56e2-7; Zeyl's slightly modified).

In a manner akin to a battle, the primary effort is concentrated in the front ranks. [17] However, this does not negate the importance of depth. Just as it is easier to overcome a single row of hoplites than many, a large cluster of the same type of particles 'attacking' from all directions, will readily 'defeat' a smaller cluster of particles of another type. Consequently, the latter would ultimately be crushed into their basic triangles (57a7-b7).

### 2. A useless movement?

---

*to another, and sometimes resting on one axis and sometimes, by revolving, on several axes.*' (X, 893d6-e1; transl. R. Bury, Harvard Univ. Press, 1967-8).

[16] For the tactics of Greek armies in Antiquity, see for instance, Hanson (1989), chapter 12.

[17] A famous example is the Battle of Thermopylae, where a small Greek army was able to repel the vastly superior Persian forces in a frontal fight. There is a wealth of texts about the Battle of Thermopylae, with the main source being Thucydides' *History*; see for instance Cawkwell (2005).

According to Timaeus' account, all parts of the universe move with the circle of the Same.[18] However, when bodies move together, they appear motionless to each other. This was well-known in Antiquity. For instance, in some ancient Greek cosmologies, such as that of the Pythagoreans, the Earth moves around a central fire and rotates on its axis.[19] Yet, since we move with it, we do not perceive its motion. A few years after Plato, one of his disciples, Heraclides, proposed a model of the universe in which the 'sphere of the fixed stars' is at rest while the Earth rotates.[20] Once again, the Earth appears absolutely still to observers moving with it. More generally, the relativity of the perception of motion was understood, as exemplified by the case of two ships at sea: to an observer on one ship, it is impossible to directly discern which ship, if any, is in motion. Hence, if all the parts of the universe are moving together, they would appear immobile to one another. Since there is nothing outside of the universe, such movement seems pointless, as summarized by Taylor:

> '*If the condition of things were absolutely uniform every state of the world should be exactly like every other. Hence all things would appear to be at rest. You might, to be sure, speculate on the possibility that the whole scheme of things was moving uniformly, but as there would be no means of distinguishing such uniformity from absolute rest, the speculation would be unprofitable.*' (Taylor (1928), p. 396).

It would be much simpler to have only the 'sphere of the fixed stars' rotating independently from the other bodies in the universe, as in later cosmologies. Why, then, does Timaeus introduce such an unnecessary complication? As we will argue in the next section, the motion of the universe, far from being useless, plays a fundamental role in Timaeus' cosmology. Indeed, while it is impossible to show directly that such a movement exists, any rotation entails a force, as is made clear for instance by the stretched rope of a rotating sling. In fact, it is due to this motion that the universe may be entirely filled by particles without any void. Because of the force resulting from the rotation, the particles are continuously pressed against each other, leaving no void:

> '*Once the circumference of the universe has comprehended the <four> kinds, then, because it is round and has a natural tendency to gather in upon itself, it constricts them all and allows no empty space to be left over.*' (58a4-7).

## IV. The rotation of the universe
### 1. The force of rotation

Since the entire universe is rotating around its polar axis, all particles within experience a rotational force. This force increases quickly with the speed of rotation, similar to the effect observed with slings.[21] Aristotle notes that, 'according to the mathematicians' of his time, the universe is vastly larger compared to the earth, which is four hundred thousand stadia. These conclusions were based on observations of stars from different locations on Earth (*On the*

---

[18] 40b1-2; 36c-e; 37c6-7; 39a1-2; 34a1-4, b3-5. Cf. also the explanations, among others, of Cornford (Cornford (1937), p. 75-76) and Taylor (Taylor (1928), p. 76). Many modern scholars, following Martin's claim (Martin (1841), p. 88), consider the Earth exempt from this movement (Dicks (1970), p. 136), while others, such as Cornford, argue that nothing in the *Timaeus* suggests such an exceptionality (Cornford (1937), p. 120-134). In any case, this does not significantly affect the overall argument, as the exception would be limited to a tiny part of the universe. Moreover, Timaeus' argument is not based on the force of rotation but rather on the *pressure* induced by it. This pressure is transmitted throughout the entire universe by the particles filling it, similar to how a ball is pressed from the outside.
[19] Aristotle, *On the Heavens*, 293a17-27; it is often assigned to Philolaus based on a fragment of Theophrastus (DK 44A16); cf. also Huffman (2007), p. 57; Huffman (2020), §4.1.
[20] For instance, Eastwood (1992).
[21] It is proportional to the square of the speed — an anachronistic result, of course.

*Heavens*, 297b30-298a20). However, Aristotle does not provide a specific measurement, and the first known measurement is found in Archimedes' texts. [22]

Although this calculation was done a century after Plato, it offers insight into how large ancient astronomers believe the universe to be. In Plato's *Statesman*, the Eleatic 'Visitor' (or 'Foreigner', 'Ξένος') describes the universe as 'immensely large' ('τὸ μέγιστον ὂν', 270a6-8). The universe completes a full rotation in a single day and night, [23] and due to its immense size, its rotational speed (cf. *supra*, note 28) is extremely fast near the 'equator' and decreases as one approaches the polar axis. Consequently, this results in significant pressure on all particles throughout the universe. [24]

When a body moves along a circular trajectory, a force acts at every point of that trajectory. Modern commentators have debated whether this force is centripetal (directed toward the center of the trajectory) or centrifugal (directed to the opposite direction). Some argue that, depending on the answer, the motions of the components of the universe would differ. [25] While this is an extremely interesting question, the debate is rooted in modern physics interpretations. Therefore, we will not discuss it, but rather focus on how ancient Greeks considered this issue.

## 2. The problem in the *Timaeus*

Circular forces have been understood since very early times, as seen in the use of slings, a tool for fighting and hunting. As the string tightens with the speed of rotation, the tension of the cord is directed away from the hand holding the sling, which is the center of the circular movement. The faster the rotation, the greater the tension in the cord. This is also evident because when the speed of rotation increases, the stone travels further and strikes its target with greater force. The same phenomenon occurs when a chariot driver takes a sharp turn, whether in battle or in a race, something apparent to anyone who has ever been in such a vehicle. [26] Anaxagoras used the same model to explain the separation of different kinds of things, emphasizing that the extreme speed of rotation creates an immense force on all rotating bodies:

> *... as these things rotated thus and were separated off by the force and speed (of their rotation). And the speed creates the force. Their speed is like the speed of nothing that now exists among men, but it is altogether many times as fast.* [27] (Frag. 9, Simplicius. *Phys.* 9, 35, 14 in Kirk-Raven (1957), 505, p. 373-4*).*

In particular, this very force is responsible for the separation of the stars, the sun, the moon, and all other bodies in the world:

> *'And all things that were to be, all things that were but are not now, all things that are now or that shall be, Mind arranged them all, including this rotation in which are now rotating the stars, the sun and moon, the air and the aither that are being separated off. And this rotation caused the separating off. And the dense is separated off from the rare, the hot from the cold, the bright from the dark and the dry from the moist.'* (Fr. 12, Simplicius *Phys.* 164, 24 and 150, 13; in Kirk-Raven (1957), 503, p. 373)

In modern terms, the ancient Greeks believed that rotation induces a centrifugal force on rotating bodies, and the faster the rotation, the greater the force.

---

[22] *The Sand Reckoner* (cf. Heath (1897), p. 226-7)
[23] That is 24 hours in modern terms.
[24] 58a; cf. also, Cornford (1937), p. 242-3 and the note of Lamb (1925) on this passage. For a more recent analysis, see Gregory (2001), p. 217-8.
[25] For instance, Archer-Hind (1888), note on 57d-58c, p. 207-8 and note 16, p. 209; Mondolfo (1934), p. 315; Cornford (1937), p. 243, note 1 and p. 244.
[26] Cf. for instance, Cornford (1937), p. 242-247; Taylor (1928), note on 57d-57d7, p. 397.
[27] 'οὕτω τούτων περιχωρούντων τε καὶ ἀποκρινομένων ὑπὸ βίης τε καὶ ταχυτῆτος. βίην δὲ ἡ ταχυτὴς ποιεῖ. ἡ δὲ ταχυτὴς αὐτῶν οὐδενὶ ἔοικε χρήματι τὴν ταχυτῆτα τῶν νῦν ἐόντων χρημάτων ἐν ἀνθρώποις, ἀλλὰ πάντως πολλαπλασίως ταχύ ἐστι' (Frag. 9, Simplicius. in *Phys.* 9, 35, 14 in Kirk-Raven (1957), 505, p. 373-4).

### 3. Forces of rotation and non-uniformity

The existence of centrifugal force (along with the debate between supporters of such a force and the those advocating for centripetal force) in a sphere rotating around its axis is inferred from the case of circular rotation around a center, as seen with a sling. However, while the rotation of a sphere does generate a force, it is neither centrifugal nor centripetal. The different parts of a sphere rotating around its polar axis follow trajectories along circles orthogonal to the axis. If we consider the axis as vertical, these circles are parallel to the 'equator' of the sphere, as shown below.

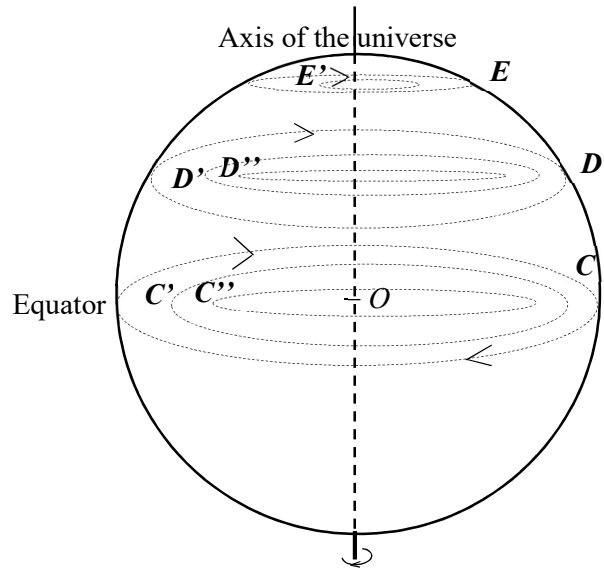

Figure 2

The points on the circles *C, C', C'', D, D', D'', E, E',* … (above figure) move along these circles, as the sphere, centered at *O,* rotates around its axis. Thus, the forces are directed away from the centers of these circles — *C, C', C'', D, D', D'', E, E',* and so on — not away from the center of the sphere. Since all points on the sphere complete one full rotation in the same amount of time, the speed along the circles increase with their radius. For instance, the speed along circle *C,* is greater than along *C', C'',* but also than along *D, D', E, E',* and so forth. Therefore, although often overlooked in commentaries, the speed [28] is not uniform throughout the universe. It is much faster near the 'equator', than near the polar axis. Accordingly, the force induced by the rotation of the sphere increases with the radius of the circular trajectory, meaning the distance from the axis of rotation. For instance, the force exerted on a particle moving along the circle *C* (the 'equator') is greater than that on a particle on *C',* which is greater than that on a particle on *C''.* The force on C'' is greater than that on a particle on *E',* and so on.

In summary, the forces due to the universe's rotation are very strong near the equator, but almost negligible near the axis. These variations in speed across different parts of the sphere play a key role in ensuring the continuity of the universe's motion (57e2-58a2).

### 4. The effect of the rotational force

Taylor and Cornford suggest that such a centrifugal force would lead to the breakdown of the universe, as its different parts will tend to escape from it. As seen in section 3 above, the forces of rotation are not centrifugal with respect of the center of the spherical universe. However, they are indeed directed away from the centers of circles of rotation, and thus away from the axis of rotation (Figure 2 above), i.e., outside of the sphere of the universe. Hence, the bodies in rotation will indeed tend to move outside of it. Yet, the claim of a breakdown resulting from such a scenario is unfounded. It is based on a false premise about the universe. As already discussed, there is nothing outside of the universe, not even space; therefore, nothing can move beyond it. [29] Hence, the rotating bodies cannot escape from the universe, which results in the forces induced by rotation exerting significant pressure on all parts of the universe. Moreover, as discussed in the previous section, since the pressure is non-uniform, the bodies, in turn, exert pressure against each other. Consequently, pressure is applied on the bodies in all directions. It is not the universe that will be broken apart, but the bodies it contains. As a result, the four regular polyhedral basic particles, namely tetrahedra, octahedra, icosahedra and cubes, are continuously destroyed and decomposed into their constituent basic triangles (56d). Consequently, due to the absence of void in the universe, new particles continuously form from these two right basic triangles, instantly replacing the destroyed ones (53c-e; 54c-d; 56d-e).

---

[28] That is, roughly speaking, the comparison between the speeds of objects is made by considering the distances covered by those objects in a given time.

However, this is not the whole picture as described by Timaeus, nor it is what we see actually around us. A universal destruction of everything by everything is not consistent with the evidence of the relative stability in the universe, nor with the goodness of the *demiourgos*. Yet the pressure on the particles is universal and affects all particles and bodies (58a-c). How, then, is it possible that the universe is not continuously destroyed and reassembled?

## V. The relative stability of the universe
### 1. The gatherings of particles

In fact, not only does the universe not go through a series of destructions and rebirths, [30] but it is everlasting (36e, 41a-b). This results from the aforementioned general rule that the same does not combat the same, much less destroy it.

> *'The reason is that a thing of any kind that is alike and uniform is incapable of effecting any change in, or being affected by, anything that is similar to it.'* (57a2-5).

This universal rule governs everything. Even without the mathematical ordering established by the *demiourgos*, the four differentiated primary bodies existed, at least in trace form (53b), these traces followed the rule that the same attracts the same, while the different distances itself from, and even combats the different (58a). Just as soldiers in an army protect fellow soldiers and combat enemies, [31] particles likewise protect those of the same kind and engage with particles of other kinds. Thus, the basic particles of a primary body can only be broken up by particles of another primary body. This can be understood physically as follows: pressure is applied to the surface of a set of particles, so that it essentially depends on their surfaces, while its resistance mostly depends of its volume. Roughly speaking, while pressure increases as the square of the number of particles, resistance increases as the cube of the number. When pressed against each other, particles of same kind, that is, of same form, will probably slide along their faces rather than crush each other. Conversely, when surrounded by particles of another kind, they will move in a disorderly manner, and the sides of the latter will cut the faces of the former, especially when the angles at the edges of the surrounding particles are small, hence sharp. In fact, they act like knife blades, easily breaking the other particles into basic right triangles (56a7-b1; d1-e2). Thus, according to the universal trend that like attracts like, particles of the same kind tend to coalesce into deep clusters, making it difficult for particles of another primary body to crush them. [32]

---

[29] Cf. Aristotle's *Physics* III, 8, 208a10-20; Alexander of Aphrodisias, *Quaestiones* III.12, 106, 35-107, 4 = DK 47 A24 = LM 14 D21); Simplicius, *On Aristotle's Physics*, 467.26-32 = DK 47 A 24; cf. also Huffman (2005), p. 549-551).

[30] Compare with the recurrent catastrophes described in the *Statesman* (270c-d), *Timaeus* (22c-23a), and *Critias* (111a-b). However, these catastrophes affect only humans, or more generally, some or even all living beings on Earth, rather than the entire universe.

[31] Cf. III.1, *supra*. In fact, Plato asserts that the establishment of the city results from a search for the protection of human beings, as it enables the city's inhabitants to assist and defend one another against both natural dangers and other people (*Republic* II, 369b-c; *Laws* III, 680d).

More generally, all of the four primary bodies, fire, air, water and earth, coalesce in different parts of the universe. Each primary body has its 'own' or 'natural' region where its particles trend to move:

> *they <the particles of each primary body> all shift, up and down, into their own proper regions* ('πρὸς τοὺς ἑαυτῶν τόπους', 58b8).

This is already true, at least partially, without the intervention of the *demiourgos* because the general disorderly shaking (53e-54b), so that the different kinds were agitated, which

> *explains how these different kinds came to occupy different regions of space, even before the universe was set in order and constituted from them at its coming to be.* (53a6-7)

In the geometrically ordered universe, the agitation continues resulting again in a separation of the four different kinds in four different regions (57c2-6), but in a much more regular manner.

### 2. The four regions of the universe

Timaeus asserts that the separation of the four natural regions depends on the mobility and stability of the four primary bodies (55d-56b). The most stable ('στάσιμος', also 'heavy', 'slow'), which are also the most difficult to move, are the cube-shaped particles of earth, because their right angles provide them with a 'natural stability' (55d8-55e7). Among the other three primary bodies, the least 'stable', which are also the most 'mobile' ('εὐκίνητος', also 'light', 'quick') particles are those with the fewest faces and thus the sharpest angles: the tetrahedral particles of fire. The particles of air and water fall between the most and least mobile, in this order, between the particles of fire and of earth (56a6-56b6).

This implies that the geometry of the four regions arises from the rotational effect of the universe on particles of varying mobility. Because the particles of fire are the swiftest ('εὐκινητότατον', 56a7), their speeds are the highest, positioning them as close as possible to the 'Equator', where the speed is greatest. Conversely, the particles of earth are positioned as close as possible to the axis of rotation, where the speed is lowest (cf., 3 *supra*, in particular Figure 2). Interposed between them are the particles of the two other primary bodies, air and water. This conclusion regarding the effect of rotational forces appears to have been widely accepted among ancient Greek thinkers, as evidenced by various fragments of Anaxagoras that emphasize the separating power of rotating forces within a 'whirl'. [33] Aristotle, in his review of the doctrine of his predecessors regarding the immobility of the Earth at the center of the universe, particularly Anaxagoras, Anaximander and Democritus, notes that they all reasoned based on the model of the 'whirl':

> *If, then, it is by constraint that the earth now keeps its place, the so-called 'whirling' movement by which its parts came together at the center was also constrained. The form of causation supposed they all borrow from observation of liquids and of air, in which the larger and heavier bodies always move to the center of*

---

[32] This was also the view of Democritus, according to Sextus Empiricus, who links it to Plato's views in the *Timaeus* (Democritus DK 68B164, Bett transl., p. 25).
[33] B12, 15-25 quoted in Simplicius' commentary on Aristotle's *Physics* 156.13-157.4; B15 *ibid.* 179.3-6; B16, *ibid.* 179.8-10 and 155.21-23. Cf. also the commentaries on these fragments in Sider (2005), note 12, p. 133-4, note 2, p. 149, note 4, p. 153-4.

*the whirl. This is thought by all those who try to generate the heavens to explain why the earth came together at the centre.* (*On the Heavens* II, 13, 295a8-15).

However, while Aristotle presents this conclusion as an empirical fact observed in nature, for Timaeus, it arises from the differences in mobility and stability of the bodies. According to this criterion and our analysis in section 3 (see above), the four natural regions of the particles of the four primary bodies are represented as shown in the figure below:

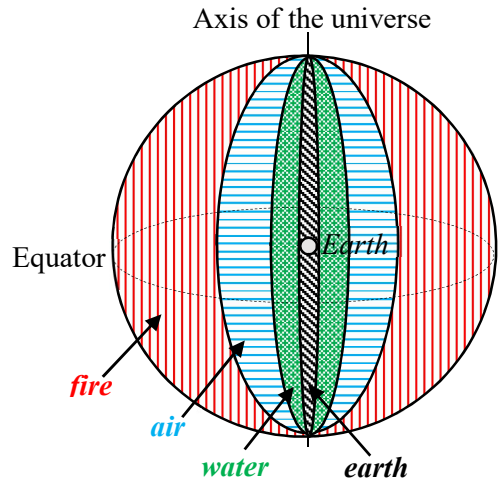

Figure 3

Note that the above representation cannot be equated with the behavior of bodies of different densities in a centrifuge, or more generally, in any device that uses rotation to separate the elements of a mixture. Timaeus does not refer to their densities, but rather their velocities. Furthermore, in a centrifuge, the lightest bodies move towards the center, while the densest gather around its walls. Accordingly, fire would be located near the center of the universe, and earth near the periphery, with air and water, in this order, in between. This is the opposite of *Timaeus'* universe, as well as almost all the representations of the universe in Antiquity. [34]

Figure 3 above is, of course, a rough scheme, as there are, for instance, other kinds of particles within the natural region of each primary body. [35] However, the main point of the figure lies in its deviation from the common representations of the four primary regions. Typically, the natural 'regions' of the four primary bodies are thought to be four concentric spherical rings of fire, air, water, and earth, arranged from the periphery to the center. [36] Cornford acknowledges that there is nothing in the text to suggest such an arrangement but argues that it can be inferred from other cosmologies (Cornford (1937), p. 246). He supports his claim with Aristotle's statement from *On the Heavens* (II, 13, 295a8-15), using the model of the whirl. According to this model, the heaviest and largest particles would be located around the center. Therefore, the regions of particles of earth would form a smaller sphere at the center, while the lightest and smallest particles would reasonably be located in the opposite direction. Consequently, the particles of fire would form a spherical ring around the circumference, centered on the universe's center, while particles of water and air would form two other spherical concentric rings in between. This is roughly how contemporary scholars view the universe as divided into four distinct regions corresponding to the four primary bodies (see, for example, Figure 3 in Gregory (2001), p. 21). However, the whirl model is inaccurate here, because a 'whirl', a concept which never appears in the *Timaeus*, is thin like a disc, [37] while the universe is not a disc but a sphere. Thus, as depicted in Figure 2 above, the speed of the particles does not depend on their distance from the universe's center, but rather on their distance from its axis of rotation.

---

[34] There are a few exceptions, such as the Pythagorean universe, where the center is occupied by a central fire.
[35] 58a-d.
[36] See Lamb (1925), 57c, note 1; Cornford (1937), p. 246; see also Miller (2003), p. 147; Aristotle's *On the Heavens*, 295b2, *ff.*

Because the Earth is close to the polar axis, its speed — and thus the centrifugal force and pressure — is much smaller than near the equator. This means that it is easier for different kinds of particles to coexist around the Earth, though not very peacefully. While the Earth is primarily composed of earth, there are also seas and rivers consisting mainly of water, as well as an atmosphere primarily composed of air. Different types of particles can even combine to form complex structures such as living beings. However, all these particles are constantly in motion, and at the interface of two large clusters of particles, they engage in conflict, each attempting to overcome the other. This dynamic is evident in phenomena such as the formation of alluvium, or, conversely, soil erosion along rivers and beaches. All these phenomena were well-known in Antiquity, as evidenced by Herodotus' account of the floods of the Nile, which he references in relation to the origins of geometry (*History*, II, 109). Such phenomena would be much less common in most parts of the 'natural' region of fire due to the significant force of rotation and the resulting high pressure, especially away from the polar axis.

### 3. A universe without void

One of the few points in their cosmologies on which both Plato and Aristotle agree is the absence of void in the universe. They even concur that this absence is related to motion. However, the reasons they provide for this are different. Aristotle argues that in a void, the speed of a body would be infinite, which he claims is absurd (*Physics*, VIII, 9-11). Thus, Aristotle's evidence for the absence of a void is a *reductio ad absurdum*. For Plato, the absence of a void results from the rotation of the universe. Due to the immense radius of the universe, which completes a full revolution in one day and night — or in modern terms, 24 hours — the pressure induced by its rotation is enormous (cf. IV.4, *supra*). As a result, all bodies are strongly constrained against one another. From the standpoint of basic particles, this implies that some are destroyed due to the pressure exerted by others, while, as soon as some spaces are freed by these destructions, new particles are formed (58a *ff*). This fundamental characteristic of *Timaeus*' universe is often overlooked. However, it is a key feature that distinguishes it from Aristotle's, and more generally, from almost all other cosmologies of Antiquity.

### 4. The mathematical paradox in *Timaeus*' universe

#### i)     Aristotle objection

Timaeus' construction of the universe encounters a significant challenge rooted in geometry. As Aristotle argues, it is mathematically impossible for several regular polyhedra to completely fill a space. [38] The meaning of 'space' in this context is not entirely clear, but it likely that Aristotle refers to the possibility of *filling* (not tiling) small spheres. The problem is whether by joining some of these polyhedra along their faces at one of their vertices, they can completely fill a sphere centered at this vertex without leaving any gaps. This construction is related to the proof of the existence of only five regular polyhedra. For easier understanding, let us consider the two-dimensional case. The question then is whether, by joining some triangles along their edges at one of their vertices, they can fill a circle without leaving any gaps.

---

[37] The model of the 'whirl' stands in contrast to the models of a cylinder or a drum, which were indeed used by other thinkers in Antiquity.
[38] *On the Heavens* III, 7, 305b27-306a20.

For instance, the answer is affirmative if we consider equilateral triangles, as six such triangles will fill the space as follows:

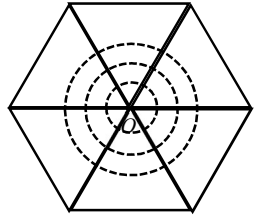

Figure 4

As seen in the figure above, the six equilateral triangles completely fill the circles of center *O* without leaving any gaps. However, they do not fill all circless centered at *O*, especially those with large radius.

Now, let us consider the isosceles triangles of angles 110°, 35°, and 35°. These triangles do not fill any circle centered on one of the vertices *O* of their largest angle:

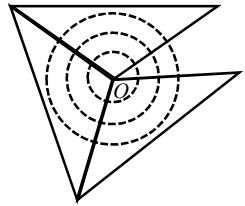

Figure 5

Indeed, three such triangles will leave some empty space in any circle centered at *O*, and it is not possible to add a fourth triangle, as there is not enough space, as shown in the figure above. The reason is clear: to fill the circle, the sum of the angles at *O* from all the triangles must be 360°. However, the sum of the angles for three such triangles is 330°, while four of them would sum to 440°, making it impossible for any number of these triangles to fill the circle. [39]

Aristotle notes that only cubes and tetrahedra can completely fill space without leaving any gaps. [40] For example, smaller cubes can do it, and can even tile a larger cube, just as squares can tile a larger square. However, this is impossible for the other regular polyhedra, [41] namely the particles of air and water. Therefore, it is impossible to fill a space with all four of Timaeus' regular polyhedra without leaving voids between them — let alone tile a sphere with these polyhedra, just as it is impossible to tile a circle with triangles or construct its boundary from straight lines. Thus, the four basic particles cannot fully occupy *Timaeus*' spherical universe, implying that there must be some voids within it. In other words, the assertions that the universe is filled with polyhedra and that it contains no voids are mathematically contradictory.

### ii) The usual "solution"

---

[39] However, if one accepts the use of similar triangles of different sizes, this is possible in the planar case because the sum of the angles of a planar triangle is equal to two right angles. In three-dimensional space, however, even if regular polyhedra of different sizes are allowed, only the cubes can 'fill space' in the sense described above, meaning they can fill small spheres centered on one vertex.

[40] *On the Heavens* III, 8, 306b8. In fact, this is not true for tetrahedra; some voids still remain between them, though they are very small.

[41] This is, however, possible by mixing some of them together, but not all of them.

As with many difficulties encountered in Plato's texts, there are typically two straightforward explanations: either Plato was unaware of the problem, or he simply did not care. However, Timaeus repeatedly asserts that there *is* no void in the universe (58a7; 59a1; 60c1; 79b1-b3, c1; 80c2) while also emphasizing the necessity of interstices ('διάκενα') (58b5; 60e5; 61a5, b1, b4) between particles. Additionally, Aristotle's criticism of the *Timaeus*, which focuses on the impossibility of the universe being both devoid of void *and* filled with polyherdra — a notion likely well-known in the Academy during Plato's time — strongly suggests that Plato was neither unaware of nor indifferent to this issue. [42] Some modern commentators, more sympathetic to Plato, have suggested an alternative hypothesis, which is not so different. Rather than attributing it to carelessness or ignorance, they suggest that this was Plato's deliberate choice, as he was not concerned with addressing minor voids (see O'Brien, (1984), p. 360-362; cf. also p. 95). According to this view, when Plato asserts the absence of voids in the universe, he is referring only to large voids. Therefore, they argue, there is no contradiction between this assertion and the existence of voids (or interstices) between the basic particles, as these voids are minuscule. This hypothesis has been accepted by most recent scholars. [43]

However, Timaeus states multiple times that there *is* no void, large or small, in the universe (58a7; 59a1; 60c1; 79b1-b3, c1; 80c2). In fact, nothing in the *Timaeus* that allows us to differentiate between voids based on their size. Moreover, this interpretation raises the question of how one would distinguish between voids that are permissible and not. At what size would a void considered small enough to be acceptable? Plato's texts offer no guidance on this matter, so any answer would be arbitrary. Additionally, the voids between particles are not always small; in fact, they are often larger than some particles (58a-b). The strongest argumen argument against this interpretation, however, comes from the fundamental premise of *Timaeus*' cosmology, which asserts the impossibility of any void whatsoever. [44] To summarize, it is impossible for the four primary particles to completely fill space; indeed, given their vast number, there are nearly infinitely many voids. Must one, therefore, concede the presence of a blatant contradiction within *Timaeus*' cosmology, as Aristotle suggests? The remainder of this article is dedicated to answering this question.

### VI. Dynamical solutions
#### 1. Necessity of the interstices

Before attempting to provide a solution to the aforementioned mathematical problem, let us first examine another issue that may seem to lead to a physical paradox. According to Timaeus' account, all the basic particles are in constant motion, driven by the movement of the universe. [45] How, then, can particles move within a universe that has no void? For instance, in a large cube completely filled with smaller cubes, the smaller cubes cannot move within the larger one. The absence of void is precisely why movement appears impossible for Parmenides' being:

---

[42] Cf. for instance Gregory (2008), note on 58b, p. 146. For a detailed criticism of such hypotheses, see Gregory (2001), p. 217

[43] For instance, see Gregory (2001), p. 217.

[44] More recent "solutions" have been proposed, such as the idea of using the soul, which is not sensible, even though it is considered a mix of the intelligible and the sensible (35a1-6), or the χώρα, which also is not sensible and is defined as a third genus, besides the intelligible and the sensible, whose 'nature is to receive all the bodies' ('περὶ τῆς τὰ πάντα δεχομένης σώματα φύσεως', 50b8) the sensible bodies. Since both seem to extend throughout the entire universe, these scholars suggest that one or the other fills these voids, addressing a physical problem through non-sensible means. In fact, neither can play any role within the universe without void space (58a7) in relation to the pressure induced on and by bodies (59a1, 60c1, 79b1-b3, 80c2), which is the issue at hand. Anyway, it is unclear how seriously the authors themselves consider their own suggestions.

[45] Cf. II.2, *supra*

> '*Nor <is there> at any spot anything more, which might keep from holding together, nor <at any one spot is there> anything less; it is instead all <of it>, full of being. Therefore, it is, all <of it>, continuous, for being keeps close to being. And so, without movement, in the bonds of great chains, it is without beginning <and> without ending, since coming into being and destruction have been driven right away, and true conviction has flung <them> afar. Staying both the same and in the same <place>, it lies by itself and stays thus fixedly on the same spot.*' (frag. B8, 23-29; O'Brien (1987)).

On the contrary, for Plato, all bodies in the universe are in motion. We will argue that, once again, the solution of this apparent paradox lies in the dynamics of the universe. Ancient Greek astronomers knew that the size of the universe was vast compared to the size of the Earth. [46] The simultaneous destruction and formation of the four primary particles that fill the universe out of and into the two basic right triangles (III.1, *supra*) result from the immense pressure generated by the rotation of the vast sphere of the universe. [47] When particles are decomposed into basic right triangles in one part of the universe, others are *simultaneously* composed from identical basic right triangles. Any possible space resulting from the decomposition of some particles is instantly filled by the composition of others. [48] These changes are governed by the so-called 'laws of transformation'. [49]

Volumes are not conserved by these 'laws'. [50] When volume of the new particles is less than that of the former ones, some void appears, which is instantaneously filled by other particles. Conversely, when the volume of the new particles is larger than the volume of the original ones, some particles are simultaneously crushed and destroyed to make room to the bigger ones. However, this could not happen if the particles completely filled the entire universe without leaving any gap, as in the case of cubes that perfectly filling a larger cube. In such a scenario, the particles would be unable to move and everything in the universe would be at rest.

In fact, Timaeus claims several times that there are many interstices between the particles, whether they are cubes or other regular polyhedra (58b5; 60e5; 61a5; b1,4). Therefore, contrary to Aristotle's claim, the existence of 'interstices' between particles is not an insurmountable difficulty arising from the mathematical impossibility of completely tiling space in Timaeus' cosmology (cf. V.4, *supra*). On the contrary, the interstices enable the movement of smaller particles between larger ones, as the size of the 'interstices' depends on the size of the particles rather than their shapes. Even if the basic particles could tile 'space' (in Aristotle's sense), there would still be some voids between them due to the irregularity of their motions. For example, while only the cubic particles of earth among the particles of the four primary bodies could tile 'space' on their own, they are the ones with the largest 'interstices':

> '*Neither air nor fire will dissolve masses of earth, because air and fire consist of parts that by nature are smaller than are the gaps ('διακένων') within earth. They thus pass without constraint through the wide spaces ('εὐρυχωρίας') of a mass of earth, leaving it intact and undissolved.*' (60e4-7).

---

[46] Cf. IV.1, *supra*.
[47] Cf. IV.4, *supra*.
[48] The model for such instantaneous changes is a rotating wheel, where the movement of one part triggers simultaneous movement in all others (cf. 79b-80c). The assertion of the simultaneity of two events — i.e., instantaneous change even when one is the cause of the other — might surprise a modern post-Galilean reader who is accustomed to considering the necessity of some interval of time between two such events (see, for instance, Gregory (2001), p. 218, on the decomposition and recomposition of atoms). However, this was not necessarily the case in Antiquity, as evidenced by discussions of phenomena as the freezing of water (Aristotle, *Sense and Sensibilia*, 417b2-3) or the instantaneous perception of light (Aristotle's *On the Soul*, 418b23), for example. This issue is explored in greater detail in a forthcoming article.
[49] These 'laws' are the mathematical rules of transformation whereby particles of fire, air, and water to change into each other, while particles of earth cannot transform into any other kind. These rules depend exclusively on the number of faces, or more precisely, on the number of half-equilateral triangles composing their faces, which must remain invariant. However, the volume of the particles changes when some particles of one kind transform into some particles of another kind. Although the total number of triangles composing the faces of the particles remains constant, the overall volume of these particles can be larger or smaller than the volume of the original ones. See for instance, Cornford (1937), p. 245-246; also Brisson-Ofman (2022).
[50] This is why many modern commentators, such as Bruins (1951), consider them a major physical flaw in Timaeus' universe.

Moreover, the fundamental law governing the universe states that 'the same moves towards the same and away from the different'. [51] This implies that particles of a given primary body are attracted to particles of the same primary body and repelled by particles of different kinds. However, this principle would lead to the complete separation of the four primary bodies. All particles of fire, air, water and earth would exist in principle would lead to distinct parts of the universe, each in their respective proper regions (cf. Figure 3, *supra*), without interacting with particles of other types:

> '*We have not explained, however, how it is that the various corpuscles have not reached the point of being thoroughly separated from each other kind by kind, so that their transformations into each other and their movement <towards their own regions> would have come to a halt.*' (58a2-58a4).

The situation, which would result in a completely motionless universe, is avoided due to the presence of the 'interstices' between basic particles, allowing particles from other primary bodies to pass through them. This is easier for the smallest particles; hence, fire particles are found almost everywhere, air particles between those of water and earth, while water particles are located between those of earth. On the other hand, it is rare to find earth particles between other particles, as they are the largest: [52]

> '*This is why fire, more than the other three, has come to infiltrate all of the others, with air in second place, since it is second in degree of subtlety, and so on for the rest.*' (58b1-2).

Due to the pressure induced by the rotation of the universe, [53] the various particles must move within the 'interstices' left by larger ones:

> '*Now this gathering, contracting process squeezes the small parts into the gaps inside the big ones.*' (58b4-5).

Conversely, when the 'interstices' are too small to allow other particles to pass through, some of the original particles are broken (cf. 4.4, *supra*). Therefore, the universe is not an inert body: the 'interstices' not only do not appear as a kind of anomaly, but they also play a fundamental role in Timaeus' cosmology.

### 2. The vanishing voids

However, while Timaeus' cosmology may not directly contradict mathematics (cf. previous section), an even more serious flaw emerges: an internal contradiction between two aspects of his account of the universe. On the one hand, Timaeus' claims the existence of 'interstices' between particles; on the other hand, he repeatedly asserts that no void *exists*. Hence the question: is it possible to reconcile both claims? Finding a solution is all the more important due to the essential role, as previously discussed, that Timaeus assigns to the interstices.

Some scholars have suggested that Plato does not concern himself with small voids, but only with large ones (cf. V.4, *supra*). However, there no hint in Timaeus' account about this claim or the significance of 'small' and 'large' voids, and this assertion is also inconsistent with

---

[51] Cf. 4.V, *supra*.
[52] This implies that the distribution of the particles of the four primary bodies is actually more complex than depicted in Figure 3, *supra*.
[53] Cf. IV.1, *supra*.

Timaeus' narrative (*Ibid.*). As noted above, the sizes of the 'interstices' between basic particles can be rather large, even larger than some basic particles themselves. For instance, the 'interstices' between particles of earth may be larger than the particles of fire or air themselves, allowing the latter to pass through them (60e4-7). Furthermore, if Timaeus did not consider the interstices due to their sizes, it would undermine their pivotal role in the order of the universe. Yet, how can this be reconciled with Timaeus' repeated assertions of the absence of void within the universe?

Based on our previous arguments, we have shown that it is improbable Plato overlooked or disregarded such an evident contradiction in the construction in *Timaeus*' cosmology (cf. V.4, *supra*). Thus, he likely believed that it was feasible for there to be no void while still allowing for the existence of nearly infinitely many 'interstices' between particles. Moreover, since he did not pose the problem or provide any explanation, he probably thought that the solution would naturally emerge from the features of *Timaeus*' universe. In fact, we believe a solution can be found along the lines of the discussion in II.2, *supra*. As the smaller particles move into the interstices left by the larger ones, either the former break the latter, or smaller particles are formed from basic triangles within these interstices, with all these particles then moving to their 'proper' regions:

> *Now this gathering, contracting process squeezes the small parts into the gaps inside the big ones. So now, as the small parts are placed among the large ones and the smaller ones tend to break up the larger ones while the larger tend to cause the smaller to coalesce, they all shift, up and down, into their own respective regions. For as each changes in quantity, it also changes the position of its region.*' (58b4-c2).

As soon as 'interstices' appear at a given location in the universe, they are *instantaneously* filled, causing their disappearance due to the movement of particles and their reconstruction from basic triangles. From a dynamic or physical perspective, the 'interstices' appear as voids — vanishing voids that are constantly filled, thus perpetually disappearing. They exist, or rather, do not exist, like some kind of ghosts, devoid of any physical existence. Within the material universe, the 'interstices' manifest as non-existent voids.

Some might argue that this is not so different from the continuous destruction and formation of particles from the basic triangles. However, there is an essential difference here. It lies in their collective behavior: as previously mentioned, a rule governs the entire universe, encompassing everything within it: the same attracts the same, while the different moves away from, or even combats, the different (cf. V.1, *supra*). This explains why the universe exists rather than breaking into basic particles (*Ibid.*). Consequently, there are four distinct regions in the universe, the four proper region corresponding to the four types of particles: fire, air, water, and earth (cf. V.2, *supra*). In contrast, the 'interstices' always remain 'interstices' — tiny, vanishing voids that never coalesce into a specific region within the universe.

This is why there are 'interstices' ('διάκενα') but no 'void' ('κενόν'); the former consistently vanish without merging into a particular area of the universe where they could then form a 'void'. Indeed, there could be no such 'void' due to the pressure induced by the universe (cf. V.4, *supra*). Hence, Timaeus' assertion that void does not exist within the universe is consistent with the presence of innumerable 'interstices' between the particles filling it. If our analysis is correct, the paradox discussed at the beginning of this article and criticized by Aristotle disappears, and *Timaeus*' cosmology is consistent regarding the problem of the void(s).

**Conclusion.** At first glance, the treatment of void in the cosmology of the *Timaeus* appears outright inconsistent. We have argued that the conventional 'solutions' to these apparent inconsistencies do not withstand careful scrutiny of the text. They portray Plato as ignorant or careless, and introduce numerous other difficulties. However, throughout this article, we have argued that these interpretations arise from viewing the universe as a static space containing all bodies, which can be either in motion or immobile. This perspective aligns with Aristotle's criticism of the *Timaeus* on the grounds that it contradicts mathematical principles.

We have attempted to demonstrate that the primary challenges posed against Plato's cosmology stem from readings that overlook the fundamental nature of its universe, as a living being in motion. In contrast to other cosmologies, *Timaeus*' entire universe is in motion. This has fundamental implications for understanding of Plato's cosmology, [54] particularly due to the significance of the universe's rotation and the different pressures it induces on various parts of the universe.

More often than not, this has been overlooked in favor of Aristotle's cosmology, leading many commentators, both ancient and modern, to analyze the *Timaeus* through Aristotelian concepts. Here, we have attempted to analyze it using Plato's concepts to demonstrate its consistency from both philosophical and physical perspectives.

## Bibliography

Unless otherwise indicated, the works of Aristotle are from his *Works* directed by W. Ross (Oxford Univ. Press). Again, unless otherwise indicated, the quotations of Plato's texts are from the translation of his works edited by J. Cooper and D. Hutchinson (Hackett Publishing, 1997), sometimes with lightly changes.

---

[54] For example, see Gregory (2001), p. 187-194, for the consequences of interpretating the treatment of the void in *Timaeus* for the receptacle (and vice-versa), although we disagree with the proposed solution.